\documentclass[11pt]{article}
\usepackage{graphicx}
\usepackage{latexsym}
\usepackage{amssymb}
\usepackage{amsmath}
\usepackage{enumerate}

\usepackage{layout}
\usepackage{eufrak}
\newtheorem{prop}{Proposition}
\newtheorem{lemma}{Lemma}

\newtheorem{corollary}{Corollary}
\newtheorem{theorem}{Theorem}
\newtheorem{remark}{Remark}

\def\real{{\mathord{{\rm I\kern-2.8pt R}}}}        
\def\inte{{\mathord{{\rm I\kern-2.8pt N}}}}

\def\sZZ{{\rm Z\kern-2.8ptem{}Z}}

\def\z{{\mathchoice
  {\sZZ}
  {\sZZ}
  {\rm Z\kern-0.30em{}Z}
  {\rm Z\kern-0.25em{}Z} }}
\def\sQQ{{\kern 0.27em \vrule height1.45ex width0.03em depth0em
          \kern-0.30em \rm Q}}
\def\qu{{\mathchoice
    {\sQQ}
    {\sQQ}
  {\kern 0.225em \vrule height1.05ex width0.025em depth0em \kern-0.25em \rm Q}
  {\kern 0.180em \vrule height0.78ex width0.020em depth0em \kern-0.20em \rm Q}
        }}
\def\sCC{{\kern 0.27em \vrule height1.45ex width0.03em depth0em
          \kern-0.30em \rm C}}
\def\complex{{\mathchoice
    {\sCC}
    {\sCC}
  {\kern 0.225em \vrule height1.05ex width0.025em depth0em \kern-0.25em \rm C}
  {\kern 0.180em \vrule height0.78ex width0.020em depth0em \kern-0.20em \rm C}
        }}

\newcommand{\R}{\mathbb{R}}

\newcommand{\ba}{\begin{array}}
\newcommand{\ea}{\end{array}}
\newcommand{\be}{\begin{equation}}
\newcommand{\ee}{\end{equation}}
\newcommand{\bea}{\begin{eqnarray}}
\newcommand{\eea}{\end{eqnarray}}
\newcommand{\beaa}{\begin{eqnarray*}}
\newcommand{\eeaa}{\end{eqnarray*}}

\def\b{\beta}

\def\z{\zeta}

\font\tenmath=msbm10 \font\sevenmath=msbm7 \font\fivemath=msbm5
\newfam\mathfam \textfont\mathfam=\tenmath
\scriptfont\mathfam=\sevenmath \scriptscriptfont\mathfam=\fivemath
\def\math{\fam\mathfam}
\def \b{\noindent}

\def \R{{\math R}}

\def\qed{ \hfill \vrule width.25cm height.25cm depth0cm\smallskip}

\newcommand{\basa}{\begin{assumption}}
\newcommand{\easa}{\end{assumption}}

\newcommand{\bas}{\begin{assum}}
\newcommand{\eas}{\end{assum}}

\newcommand{\ignore}[1]{}
\begin{document}

\renewcommand{\thefootnote}{\fnsymbol{footnote}}

\renewcommand{\thefootnote}{\fnsymbol{footnote}}

\title{Density for solutions to stochastic differential equations with unbounded drift     }

\author{Christian Olivera $^{1,}$\footnote{Supported by
FAPESP 2012/18739-0} \hskip0.2cm 
Ciprian A. Tudor $^{2}$ \vspace*{0.1in} \\
$^{1}$ Departamento de Matem\'atica, Universidade Estadual de Campinas,\\
13.081-970-Campinas-SP-Brazil. \\
colivera@ime.unicamp.br \vspace*{0.1in} \\
 $^{2}$ Laboratoire Paul Painlev\'e, Universit\'e de Lille 1\\
 F-59655 Villeneuve d'Ascq, France.\\
 \quad tudor@math.univ-lille1.fr\vspace*{0.1in}}
\maketitle

\begin{abstract}
Via a special transform and by using the techniques of the Malliavin calculus, we analyze the density of the solution to a stochastic differential equation with unbounded drift.
\end{abstract}

\medskip


 \smallskip


\section{Introduction}
Our purpose is to prove the existence and other properties for the density of the solution to the  stochastic differential equation (SDE) in $\mathbb{R}^ {d} $ 
\begin{equation}
\label{sde}
X_{t} (x)= x+ \int_{0}^ {t} b(X_{r}(x) ) dr + \int_{0} ^ {t} \sigma (X_{r} (x)) d B_{r}, \hskip0.5cm t\in [0,T], x\in \mathbb{R} ^{d}
\end{equation}
with {\it unbouded H\"older  continuous drift} $b$ and smooth diffusion coefficient $\sigma$.

The Malliavin calculus is nowadays a widely recognized mathematical   theory to  study densities of random variables in general,  and of solutions to stochastic (partial) differential equations  in particular. On the other hand, the theory works well when the coefficients of the stochastic equations are smooth enough. For instance, one can prove the absolute continuity with respect to the Lebesque measure  of the law of solutions to  SDEs   with coefficients which are globally Lipschitz continuous and of at most linear growth and with invertible diffusion matrix. For more regular coefficients, e.g. for bounded functions with bounded derivatives, we can get the existence of a smooth density for the solution to (\ref{sde}). We refer, among others, to \cite{N} or \cite{Sans} for a more complete exposition.

On the other hand, when the coefficients are less regular, for example if they  are only H\"older continuous, then much more work is needed in order to obtain the Malliavin differentiability of the solution or the existence and the smoothness of  its density. There already exists a large literature in this direction, see among others,  \cite{ba}, \cite{Ko}, \cite{Ko3}, \cite{Romi}, \cite{Ha}, \cite{Ba},\cite{Deb}, \cite{Marco} etc.  As far as we know, the study of the density for solutions to (\ref{sde}) in the case of {\it unbounded H\"older continuous drift } $b$ and in general dimension $d\geq 1$ is  not contained  in none of these references. Let us briefly discuss the closest results to our paper. For $d=1$, the absolute continuity of the solution to (\ref{sde}) has been showed in \cite{PF}. In \cite{Deb}, the existence of the density for the solution to a  special  one-dimensional SDE with H\"older continuous diffusion and drift with at most linear growth has been considered. Further improvements have been achieved in  \cite{Ba} in the general case $d\geq 1$,  by weakening the H\"older continuous hypothesis. In \cite{Marco}, the author shows that if there exists some ball in  $\mathbb{R}^{d}$  in which both coefficients  are smooth and with bounded partial derivatives, then the
density is smooth inside a smaller ball. In the case $d\geq 1$, \cite{Ha} shows that if there exists some open interval on which $b$ is H\"older continuous,  uniformly elliptic
and smooth, then the density is H\"older continuous on the interval. The case of bounded H\"older continuous diffusion coefficient in dimension $d\geq 1$ is treated also in \cite{BCC}.  Stronger, or different, assumptions on the drift are considered in  \cite{ba}, \cite{Ko}, \cite{Ko3}, \cite{Romi}.

Our purpose is to present a new and easier method to deal with densities of solutions to SDEs with irregular drift, via some special  transformations of the drift and of the diffusion coefficient. The basic ideas is as follows: by using a transformation defined in \cite{FGP}, we can show that the equation (\ref{sde}) is equivalent to another SDE
\begin{equation}
\label{2}
Y_{t}(x)=  x+ \int_{0}^ {t} \widetilde{b}(Y_{r}(x) ) dr + \int_{0} ^ {t}\widetilde{ \sigma} (Y_{r} (x)) d B_{r}
\end{equation}
where the "new coefficients" $\widetilde{b}$ and $\widetilde{\sigma}$ are  explicitly constructed  from $b$ and $\sigma$.  
The equivalence between the equations (\ref{sde}) and (\ref{2})  means that if  (\ref{sde}) has a solution $ X_{t}(x)$ then (\ref{2}) has a solution $ Y_{t}(x)$ and  $ X_{t}(x)= \Phi ^ {-1} (Y_{t} (x) ) $ for every $t\in [0,T], x \in \mathbb{R} ^ {d}$ with $ \Phi$ a smooth and invertible function. On the other hand, the new coefficients  $\widetilde{b}$ and $\widetilde{\sigma}$ are smooth and the standard Malliavin techniques can be applied to (\ref{2}). This will easily lead to the Malliavin differentiability of the solution to (\ref{sde}) and to the existence of  its probability density. This is the main purpose of our work, to show that our approach allows, by relatively trivial arguments, to obtain results that usually demand significant technical work. A similar idea has been marginally employed in \cite{Romi1}.

Our paper is organized as  follows.   Section 2  is devoted to the introduction of some notation and of the basic elements of the Malliavin calculus. In Section 3 we present the special transform used in our work and we derive some consequences for the Malliavin differentiability  and the density of the solution. In Section 4 we analyze the existence of the density for functionals of the solution.

\section{Preliminaries}
This preliminary section is devoted to the presentation of the basic tools from Malliavin calculus and to the introduction of some notation needed in the paper.

\subsection{Malliavin derivative}

\hspace{1.1cm} Let us  present the elements from the Malliavin calculus that  will be used in the paper.  We refer to \cite{N} for a more complete exposition. Consider ${\mathcal{H}}$ a real separable Hilbert space endowed with the scalar product $\langle \cdot, \cdot \rangle _{\mathcal{H}}$ and $(B (\varphi), \varphi\in{\mathcal{H}})$ an isonormal Gaussian process \index{Gaussian process} on a probability space $(\Omega, {\cal{A}}, \mathbb{P})$, that is, a centred Gaussian family of random variables such that $\mathbf{E}\left( B(\varphi) B(\psi) \right)  = \langle\varphi, \psi\rangle_{{\mathcal{H}}}$.

We denote by $D$  the Malliavin  derivative operator that acts on smooth functions of the form $F=g(B(\varphi _{1}), \ldots , B(\varphi_{n}))$ ($g$ is a smooth function with compact support and $\varphi_{i} \in {{\cal{H}}}, i=1,...,n$)
\begin{equation*}
DF=\sum_{i=1}^{n}\frac{\partial g}{\partial x_{i}}(B(\varphi _{1}), \ldots , B(\varphi_{n}))\varphi_{i}.
\end{equation*}

It can be checked that the operator $D$ is closable from $\mathcal{S}$ (the space of smooth functionals as above) into $ L^{2}(\Omega; \mathcal{H})$ and it can be extended to the space $\mathbb{D} ^{1,p}$ which is the closure of $\mathcal{S}$ with respect to the norm
\begin{equation*}
\Vert F\Vert _{1,p} ^{p} = \mathbf{E} F ^{p} + \mathbf{E} \Vert DF\Vert _{\mathcal{H}} ^{p}. 
\end{equation*}

We denote by  $ \mathbb{D} ^{k, \infty}:= \cap _{p\geq 1 } \mathbb{D} ^{k,p}$ for every $k\geq 1$. In this paper, $\mathcal{H}$ will be the standard Hilbert space $L^ {2}([0,T])$.

We will use  the chain rule for the Malliavin derivative (see Proposition 1.2.4 in \cite{N}) which says that if $\varphi: \mathbb{R}\to \mathbb{R}$ is a differentiable function with bounded derivative and $F\in \mathbb{D} ^ {1,2}$, then  $\varphi (F) \in \mathbb{D} ^ {1,2}$ and
\begin{equation}
\label{chain}
D\varphi(F)= \varphi ' (F) DF.
\end{equation}
For a random variable $F$ with values in $\mathbb{R}$, we wll say that is Malliavin differentiable if all its components are in $\mathbb{D} ^{1,2}. $

An important role of the Malliavin calculus  is that it provides criteria for the existence of the density of a random variable. This was actually the initial motivation to construct this mathematical theory.  There exist a huge number of formulas for densities of random variables in terms of the Malliavin operators, see \cite{N} or \cite{NPbook} among others.  Here we will use the following  result: 
 if $F$ is a random variable in $\mathbb{D} ^ {1,2}$ such that 
\begin{equation}\label{25o-1}
\Vert DF \Vert _{\mathcal{H}} >0 \mbox{ a.s. then $F$ admits a density with respect to the Lebesque measure }
\end{equation}
 (see e.g. Theorem 2.1.3 in \cite{N}).
This is a classical result in Malliavin calculus, related to the invertibility of the Malliavin matrix, and it goes  back to  \cite{Shige}.  More recently, this link between Malliavin calculus and the probability densities of random variables has been reinforced and new interesting results have been obtained (see, among others, \cite{Ba}, \cite{NV}, \cite{NQ}). 

\subsection{Notation}

Let us introduce some notation needed throughout the paper. 

For any  $\theta \in (0,1)$, we define the set   $C^{\theta}(\R^d;\R^k)$, $d,k \geq 1$ of   the  mappings   $f:\R^d \to \R^k$  such that 
$$
[f]_{\theta}:=\sup_{x,y\in \mathbb{R} ^ {d}, x\neq y, |x-y|\leq 1} \frac{\vert f(x)-f(y)\vert }{\vert x-y \vert ^\theta} < \infty\,.
$$
By  $\vert \cdot \vert, \langle \cdot, \cdot\rangle  $ we denote the Euclidean norm and scalar product in  $\mathbb{R} ^ {k}$ respectively.

 The space $C^{\theta}(\R^d;\R^k)$ is a  Banach space with respect to  the norm
$$
\|f\|_{\theta}=\|(1+|\cdot|)^{-1}f\|_\infty+[f]_{\theta}
$$
where $\Vert \cdot \Vert _{\infty} $ denotes the supremum norm. It is also called in \cite{FGP} as the space of "locally uniformly $\theta$-H\"older functions".

We will say that $  f \in C ^ {n+\theta } (\mathbb{R} ^ {d}; \mathbb{R} ^ {d})$ ($n\geq 1$ integer) if $ f \in C^{\theta}(\R^d;\R^d)$ and moreover the Fr\'echet derivatives $\mathcal{D} ^ {i} f (i=1,..,n)$ are bounded and H\"older continuous of order $\theta \in (0,1)$.

\section{The special transformation and the density of the solution}

Assume $k,d\geq 1$ and let $\left( B_{t}\right) _{t\in [0,T]}  $ be a standard $k$-dimensional Brownian motion ($B=(B ^ {1},..., B ^ {k})$) on a probability space $\left( \Omega, \mathcal{F}, P\right)$ with respect to its natural filtration $\left( \mathcal{F} _{t} \right) _{t\geq 0}. $ Consider the following stochastic differential equation in $\mathbb{R} ^ {d} $ 
\begin{equation}
\label{sde1}
X _{t}(x)= x + \int_{0} ^ {t} b (X_{s} (x)) ds + \int_{0} ^ {t} \sigma (X_{s}(x) ) d B_{s}, \hskip0.5cm t\in [0,T] 
\end{equation}
or equivalently 
\begin{equation*}
X_{t}(x)=x+\int_{0} ^ {t} b (X_{s} (x)) ds+\sum_{i=1}^ {k}\int_{0} ^ {t}\sigma _{i}  (X_{s}(x) ) d B^ {i}_{s}
\end{equation*}
where $x\in \mathbb{R} ^ {d} $ and the stochastic integral in (\ref{sde1}) is understood in the It\^ o sense.

As in \cite{FGP}, we shall  assume throughout  our work that the vector fields $b$ and $\sigma $ satisfy the following conditions
\begin{equation}\label{con1}
    b\in C^{\theta}(\R^d; \R^d)\,,
\end{equation}
and for every $i=1,.., k$,

\begin{equation}\label{con2}
    \sigma _{i} \in C_{b}^{3}(\R^d; \R^{d})\,, 
\end{equation}

\noindent  and $a=\sigma \sigma^{\star}$ ($\sigma ^ {\star}$ denotes the adjoint matrix of $\sigma$) is invertible  and satisfies

\begin{equation}\label{con3}
    \| a\|_{0}= \sup_{x\in \mathbb{R}^{d}} \| a^{-1}(x)\|< \infty 
\end{equation}
where for every $x\in \mathbb{R} ^ {d}$,  $\Vert a^ {-1}(x) \Vert $ denotes the Hilbert-Schmidt norm of the $d\times d$  matrix $a ^ {-1}(x) $.

\subsection{The auxiliary elliptic system}

Fix $\lambda >0$ and fix $0< \alpha < \theta $ (recall that $\theta $ is fixed by  condition (\ref{con1})). We introduce the auxiliary elliptic systems

\begin{equation}\label{5s-1}
\lambda \psi_\lambda - L \psi_\lambda= b,
\end{equation}

  \noindent where $L$  is the Kolmogorov operator associated to (\ref{sde1}) defined by

\begin{equation}
\label{5s-2}
	L= \frac{1}{2} Tr( \sigma\sigma^{\star} \mathcal{D}^{2}) + b(x)\mathcal{D}  .
\end{equation}
	
\b We also introduce the function

\begin{equation}
\label{5s-3}
\phi_{\lambda}(x)= x + \psi_\lambda(x)
\end{equation}
	where $ \psi _{\lambda }$ is the unique classical solution in $ C ^ {2+ \alpha }(\mathbb{R} ^{d}; \mathbb{R} ^ {d})$ of the equation (\ref{5s-1}) (which exists due to Theorem 5 in \cite{FGP}).

\vskip0.2cm

The following lemmas are due to \cite{FGP} (see Theorem 5 and Lemma 8 in this reference) and they are key points  for our method. Notice that we need to have $\sigma$ of class $C^{3}$ in order to have the second point below.

\begin{lemma}\label{19o-1}For  $\lambda >0$, let  $\phi_{\lambda}$ be given by (\ref{5s-3}). Then the function  $ \phi _{\lambda }$ has  the following properties when $\lambda $ is large enough:
	
	\begin{enumerate}
	\item There exists $C>0 $ such that $\| \mathcal{D} 	\phi_{\lambda}\|_0 < C.$ 
	\item  $	\phi_{\lambda}$ is a $C^{2}$-diffeomorphism. 
	\item    $\phi_{\lambda}^{-1}$  has bounded first and second derivatives.
	\end{enumerate}

\end{lemma}

In the rest of the paper, we will consider $\lambda $ large enough in order that the point 1.-3. in Lemma \ref{19o-1} abovehold.  The following special transform allows to remove the irregular drift in (\ref{sde1}). 

\begin{lemma}\label{25o-2}
We define the "new' coefficients $\widetilde{b} $ and $\widetilde{\sigma}$ from $b$ and $\sigma$ by  

\begin{equation}
\label{news}
\tilde{b}(x)= \lambda  \psi_{\lambda}(\phi_{\lambda}^{-1}(x))   \mbox{ and }  \tilde{\sigma}(x)= \mathcal{D}\phi  _{\lambda}(\phi_{\lambda}^{-1}(x)) \sigma(\phi_{\lambda}^{-1}(x))
\end{equation}
for every $x\in \mathbb{R} ^ {d}$. 
Then
\begin{enumerate}
\item $\widetilde{b} \in C ^ {1+\alpha} (\mathbb{R} ^ {d}; \mathbb{R} ^ {d})$ and $\widetilde{\sigma}\in C ^ {1+\alpha} (\mathbb{R} ^ {d}; \mathbb{R} ^ {d\times k})$.


\item Consider the SDE
\begin{equation}
\label{sde2}
Y_{t}(x) = y + \int_{0} ^ {t} \widetilde{b} (Y_{s}(x)) ds + \int_{0} ^ {t} \widetilde{\sigma} (Y_{s}(x))dB_{s}, \hskip0.5cm t\in [0,T].
\end{equation}
Then the SDE (\ref{sde2}) is equivalent to (\ref{sde1}) in the following sense: if $ X  $ solves (\ref{sde1}), then $Y$ defined by 
\begin{equation}\label{xy}
Y_{t}(x)= \phi_{\lambda } \left( X_{t}(x)) \right)  
\end{equation}
solves (\ref{sde2}) with $y=\phi_{\lambda }  (x)$. The converse implication is also true.
\end{enumerate}
\end{lemma}


The above two lemmas easily lead  to the following result. Recall that a random variable in $\mathbb{R}^{d}$ is Malliavin differentiable if all its components are in $\mathbb{D} ^{1,2}. $

\begin{theorem}\label{t2}
Assume  (\ref{con1}),  (\ref{con2}) and   (\ref{con3}).
Then, for every $t\in [0,T]$ and for every $x\in \mathbb{R}^{d}$, we have :

 \begin{enumerate}
\item  The random variable $X_{t}(x)$ is Malliavin differentiable.

 \item The random variable $X_t(x)$ 
admits a density  $\rho_{X_{t}(x)}$  with respect to the Lebesgue measure. 
\end{enumerate}
\end{theorem}
{\bf Proof: }
 Recall that $Y$ solves the SDE   (\ref{sde2}) which is a SDE with smooth coefficients since 
$\widetilde{b}, \widetilde{\sigma} \in C ^ {1+\alpha }$, see Lemma \ref{25o-2}. By standard arguments based on  Malliavin calculus (see Chapter 2 in  \cite{N}) it follows  that for every $t,x$, the random variable $Y_{t}(x)$ is Malliavin differentiable. Since $\phi _{\lambda } ^ {-1}$  is differentiable with bounded derivative  by Lemma \ref{19o-1}, we obtain the Malliavin differentiability of $X_{t}(x)$ for every $t,x$. Moreover,  by applying the chain rule (\ref{chain}),  for every $r\in (0, T]$

\begin{equation}\label{9s-1}
D_r X_t (x) = \mathcal{D}\phi_{\lambda}^{-1}(Y_t) D_r Y_t (x).
\end{equation}

To show point 2., we note that $Y_{t}(x)$ admits a density with respect to the Lebesque measure. This can be argued by using Theorem 2.3.1 in \cite{N}, based on the fact that the diffusion matrix of $Y$ is non-singular, see Lemma \ref{ll4} below.  Denote by $ \rho_{Y_{t}(x)}$ the density of the random variable $ Y_{t}(x)$.

Next, for every continuous bounded function $\varphi : \mathbb{R} ^ {d} \to \mathbb{R}$, we can write, since $\phi_{\lambda}$ is a $C^ {2}$-diffeomorphism, 

\begin{eqnarray*}
\mathbf{E} \varphi ( X_{t}(x)) &=& \mathbf{E} \varphi \left( \phi_{\lambda} ^ {-1} (Y_{t}(x))\right) \\
&=& \int_{\mathbb{R} ^ {d}}\varphi \left( \phi_{\lambda}^ {-1} (u)\right) \rho _{Y_{t}(x)}(u) du\\
&=&\int_{\mathbb{R} ^ {d}} \varphi (u) \left(\det  J \phi_{\lambda }(u) \right) \rho_{Y_{t}(x)}(\phi_{\lambda } (u))du
\end{eqnarray*}
where $J$ denotes the Jacobian (with respect to $x$). It  follows from the last equality that for every $t>0$ and for every $x\in \mathbb{R} ^ {d}$  the random variable  $X_{t}(x)$ admits a  density $ \rho_{X_{t}(x) }$ with  respect to Lebesgue and the following formula holds true

\begin{equation}
\rho_{X_{t}(x)} = \left(\det J\phi_{\lambda } \right)\cdot  \rho_{Y_{t}(x)} (\phi_{\lambda}).
\end{equation} \qed

\begin{remark}
Let us point out that in general (see e.g. \cite{Ko2}, \cite{Ko3}, \cite{Zhang3}), proving the Malliavin differentiability of solutions to SDE with irregular drift demands pretty technical work. 
\end{remark}

It is also possible to give, in dimension 1, a representation of the density in the spirit of \cite{NV}. The following result is a consequence of Theorem 3.1 in \cite{NV}.

\begin{corollary}

Assume $d=1$ and let the assumptions in Theorem \ref{t2}.  Then, for every $t\in [0,T]$ and $x\in \mathbb{R} ^{d}$, the density of  the random variable $X_t(x)$  
can be expressed as

\[
\rho_{X_{t}(x)}(z)= \frac{\mathbf{E}|X_t (x)|}{2g _{X_{t}(x)} (z)} e^{-\int_{0}^{z} \frac{u}{g_{X_{t}(x)}(u)} \ du}.
\]

\end{corollary}

\section{Density of functions of the solution}

In this paragraph our aim is to prove the existence of densities for functions of the solution to  SDEs with unbounded drift. Although the main idea comes from the special transformation (\ref{xy}), other auxiliary results   are also neeeded.

We will first prove some auxiliary lemmas.

\begin{lemma}\label{ll3}

Let $\phi_{\lambda} $ be given by (\ref{5s-3}). Then for every $x \in \mathbb{R} ^ {d}$ and for $\lambda $ large enough
\begin{equation}
\label{19o-3}
\left| \mathcal{D} \phi_{\lambda } (x) \right| \geq C >0.
\end{equation}
\end{lemma}
{\bf Proof: }From (\ref{5s-3}) we have for every $x\in \mathbb{R}^ {d}$ 
\begin{equation}
\label{19o-4}
\mathcal{D} \phi_{\lambda}(x)= I_{d}+ \mathcal{D} \psi _{\lambda } (x)
\end{equation}
where $I_{d}$ denoted the $d$-dimensional identity matrix.  On the other hand, it follows from formula (18) in \cite{FGP} that for every $x\in \mathbb{R}^ {d}$
\begin{equation}
\label{19o-5}
\left| \mathcal{D} \psi_{\lambda }(x) \right| \leq c(\lambda)
\end{equation}
where $c(\lambda)\to 0$ as $\lambda \to \infty. $ From (\ref{19o-4}), (\ref{19o-5}) and the triangle inequality, we get

\begin{eqnarray*}
\left| \mathcal{D} \phi_{\lambda } (x) \right| &=& \left| I_{d}+ \mathcal{D} \psi _{\lambda } (x)\right| \\
&\geq & 1- \left| \mathcal{D} \psi_{\lambda }(x)  \right|  \geq C >0
\end{eqnarray*}
for $\lambda $ large enough.  \qed

\begin{lemma}\label{ll4} Let $\tilde{\sigma}$ be given  by (\ref{news}). Then the diffusion matrix of the SDE (\ref{sde2})  $B= \widetilde{\sigma} \widetilde{\sigma } ^ {\ast} $ is uniformly elliptic. 

\end{lemma}
{\bf Proof: } From the definition of the coefficient $\widetilde{\sigma}$ (\ref{news}) we can write (recall that $a=\sigma \sigma ^ {\ast}$) for every $x\in \mathbb{R}  ^ {d}$ 
 
\begin{equation*}
B(x)= (\mathcal{D} \phi_{\lambda} )(\phi_{\lambda } ^ {-1})a (\phi_{\lambda } ^ {-1})\left( (\mathcal{D} \phi_{\lambda} )(\phi_{\lambda } ^ {-1})\right) ^ {\ast}
\end{equation*}
and this implies that the inverse matrix of $B$ satisfies
\begin{equation*}
B^ {-1} (x)=\left[ \left( (\mathcal{D} \phi_{\lambda} )(\phi_{\lambda } ^ {-1})\right) ^ {\ast}\right] ^ {-1} a ^ {-1} (\phi_{\lambda } ^ {-1})\left(  (\mathcal{D} \phi_{\lambda} )(\phi_{\lambda } ^ {-1})\right) ^ {-1}.
\end{equation*}
Notice that, since $a^ {-1} $ is bounded from (\ref{con3}) and by using the properties of $\phi_{\lambda}$ in Lemma \ref{19o-1}, we immediately obtain the existence of a strictly positive constant $C$ such that 
$$\vert B^ {-1} (x) \vert \leq C $$
for every $x$.

Denote by $\Vert  A\Vert $ the operator norm of the matrix $A$ and recall that if  $A=A ^ {\ast} $ then 
$$\Vert A \Vert =\sup_{\lambda \in Spec (A) } \vert \lambda  \vert $$
where $Spec (A)$ is the spectrum of $A$.  The above formula implies that
$$ \lambda _{i,B^ {-1} } (x) \leq C $$  for every $i\in I$ where $Spec(B^ {-1})=\{ \lambda _{i, B ^ {-1}}, i\in I\}.$ Consequently,
\begin{equation}
\label{20o-1}
\lambda _{j, B} (x) \geq C 
\end{equation} for every $j\in J$ where we denoted  the spectrum of $B$ by $Spec(B)= \{ \lambda_{j,B}, j\in J\}$. 
The above inequality (\ref{20o-1}) clearly gives
$$ \langle  B(x)\xi, \xi \rangle \geq C \vert \xi \vert ^ {2} \mbox{ for every } \xi  \in \mathbb{R} ^ {d}.$$

The converse inequality follows easily since $B$ is unifomly bounded and then 
$$  \langle  B(x)\xi, \xi \rangle \leq \vert B(x)\xi \vert  \vert \xi \vert  \leq C \vert \xi \vert ^ {2}.$$  \qed

\begin{lemma}\label{ll5} Let $Y$ be the solution to (\ref{sde2}). Then for every $t\in [0, T], x\in \mathbb{R} ^ {d}$ 
\begin{equation*}
JY _{t}(x)=e ^ {\int_{0} ^ {t}   (\mathcal{D} \widetilde{\sigma} ) (Y_{s}(x))dB_{s} + \int_{0} ^ {t} \left[ (\mathcal{D} \widetilde{b} ) (Y_{s}(x))-\frac{1}{2}(\mathcal{D} \widetilde{\sigma} )^ {2} (Y_{s}(x)) \right]}
\end{equation*}
and
\begin{equation*}
JY_{t} ^ {-1} (x)=e ^ {-\int_{0} ^ {t}   (\mathcal{D} \widetilde{\sigma} ) (Y_{s}(x))dB_{s} - \int_{0} ^ {t} \left[ (\mathcal{D} \widetilde{b} ) (Y_{s}(x))-\frac{1}{2}(\mathcal{D} \widetilde{\sigma} )^ {2} (Y_{s}(x)) \right]}
\end{equation*}
\end{lemma}
{Proof: } Notice that the Jacobian of $Y$  satisfies the equation

\begin{equation*}
JY_{t}(x)= 1+ \int_{0} ^ {t} (\mathcal{D} \widetilde{b} ) (Y_{s}(x)) Y_{s}(x) ds+ \int_{0} ^ {t} (\mathcal{D} \widetilde{\sigma} ) (Y_{s}(x)) Y_{s}(x)dB_{s}
\end{equation*}
while  its inverse solves (see \cite{N}, page 126, formula (2.58)), for $r\in [0,T], x\in \mathbb{R} ^ {d}$,
\begin{equation*}
J Y_{r} ^ {-1}(x)= 1- \int_{0}^ {r} (\mathcal{D} \widetilde{\sigma} ) (Y_{s}(x)) (JY _{s} ^ {-1}(x) )dB_{s} - \int_{0}^ {t} 
\left[ \mathcal{D} \widetilde{b} ) (Y_{s}(x))-( \mathcal{D} \widetilde{\sigma } ) ^ {2} (Y_{s}(x))\right] (Y_{s}^ {-1}(x) ) ds.
\end{equation*}
By solving the above two equations, we get the conclusion. \qed

For simplicity, we will assume in the sequel  that   $d=k=1$.   Consider a function $ G :[0,T] \times \mathbb{R} \to \mathbb{R} $ such that $G(t, \cdot ) \in C_{b} ^ {1}(\mathbb{R})$ and 
\begin{equation}
\label{19o-2}
\vert \partial _{x} G (t,x)\vert  \geq C>0 
\end{equation}
for every $t\in [0,T] $ and for every $x\in \mathbb{R}$, where $\partial_{x}$ denote the derivative with respect to the variable $x$. 

The main result of this section states as follows.

\begin{prop}\label{t3}
We  assume  (\ref{con1}),  (\ref{con2}) and   (\ref{con3}).
Then, for every $t\in [0,T]$ and for every $x\in \mathbb{R}$, we have :

 \begin{enumerate}
\item  The random variable $G(t,X_t(x))$  is Malliavin differentiable.

 \item The  random variable $G(t,X_t(x))$ 
admits a density  $\rho_{t,x}$  with respect to the Lebesgue measure.

 \end{enumerate}
\end{prop}
{\bf Proof: }  Fix $t\in [0, T], x\in \mathbb{R}$. We have seen in Proposition \ref{t2} that the random variable $X_{t}(x) $ belongs to $ \mathbb{D} ^ {1,2}$. Since $G$ is of class $ C_{b} ^ {1}$ with respect to $x$, we get that $G(t,X_t(x)) \in \mathbb{D} ^ {1,2}$.

Moreover, for every $r>0$,  
\begin{eqnarray*}
D_{r} G(t, X_{t}(x)) &=& \partial _{x}G (t, X_{t}(x) ) D_{r} X_{t}(x)\\
&=&  \partial _{x}G (t, X_{t}(x) ) \mathcal{D} \phi_{\lambda }^ {-1} (Y_{t}(x) ) D _{r} (Y_{t}(x))
\end{eqnarray*}
where we used relation (\ref{xy}) and  we used twice the chain rule (\ref{chain})  for the Malliavin derivative.

On the other hand, from \cite{N}, formula (2.59) on page 126, we have
\begin{equation*}
D_{r} Y_{t}(x) =J Y_{t}(x) \left( J Y_{r} (x) \right) ^ {-1} \widetilde{\sigma} (Y_{r}(x))
\end{equation*}
where $ (JY_{t}(x)) ^ {-1}$ represents the inverse of $ JY_{t}(x)$ (recall that $Y$ generates a $ C^ {1}$ flow of diffeomorphism, see \cite{FGP}). We thus get for $r\in (0,T]$
\begin{equation}
\label{20o-3}
D_{r} G(t, X_{t}(x)) =\partial _{x}G (t, X_{t}(x) ) \mathcal{D} \phi_{\lambda }^ {-1} (Y_{t}(x) ) J Y_{t}(x) \left( J Y_{r} (x) \right) ^ {-1} \widetilde{\sigma} (Y_{r}(t)).
\end{equation}
and
\begin{eqnarray*}
\Vert D_{r} G(t, X_{t}(x))\Vert ^ {2} _{L ^ {2} ([0,T] )} &=& \int_{0} ^ {T} \left| \partial _{x}G (t, X_{t}(x) ) \mathcal{D} \phi_{\lambda }^ {-1} (Y_{t}(x) ) J Y_{t}(x) \left( J Y_{r} (x) \right) ^ {-1} \widetilde{\sigma} (Y_{r}(t))\right| ^ {2}dr\\
&\geq & C \vert J Y_{t}(x) \vert ^ {2}  \int_{0} ^ {T} \left|   \left( J Y_{r} (x) \right) ^ {-1}\right| ^ {2} dr
\end{eqnarray*}
where we used condition (\ref{19o-2}), Lemma \ref{ll3} and Lemma \ref{ll4} to bound from below $ \partial _{x}G (t, X_{t}(x) )$, $
 \mathcal{D} \phi_{\lambda }^ {-1} (Y_{t}(x) )$ and $\widetilde{\sigma} (Y_{r}(t))$ respectively.
The result follows easily from Lemma \ref{ll5} and the criterion (\ref{25o-1}). \qed

\begin{remark}
If $d\geq 1$, and the partial derivatives of $G:[0, T]\times \mathbb{R} ^{d}\to \mathbb{R}$ with respect to $x\in \mathbb{R} ^{d}$  are  bounded (above and below)  by strictly positive constants, then the above result can be also obtained. The formula (\ref{20o-3}) will become 
\begin{eqnarray*}
D_{r}G(t, X_{t}(x))&=& \sum_{i=1}^{d} \frac{\partial G}{\partial x_{i}} (t, X_{t}(x)) \frac{\partial}{\partial x_{j} } (\phi_{\lambda } ^{-1})^{i}(Y_{t}(x)) \left( JY_{t}(x)\right) _{j,l}\left( \left( JY_{r}(x)\right) ^{-1}\right) _{lk} \widetilde{\sigma}^{k} (Y_{r}(x))
\end{eqnarray*}
where the superscript means the components of functions and the subscripts are the elements of matrices. By using the boundedness of the partial derivatives of $G$ and Lemmas \ref{ll3}, \ref{lll4} and \ref{ll5}, the strict positivity of $\Vert D_{r} G(t, X_{t}(x))\Vert ^ {2} _{L ^ {2} ([0,T] )} $ can be proven.
\end{remark}

\vskip0.3cm
{\bf Acknowledgement: } The authors acknowledge partial support from the CNRS-FAPESP grant 267378.
C. Olivera  is partially supported by FAPESP 	by the grants 2017/17670-0. and 2015/07278-0.

\end{document}